\documentclass[a4paper,12pt]{amsart}
\usepackage{amssymb}
\usepackage{graphicx}
\usepackage{ifthen}
\nonstopmode \numberwithin{equation}{section}
\newtheorem{thm}{Theorem}
\newtheorem{cor}{Corollary}
\newtheorem{lem}{Lemma}

\newtheorem{conj}{Conjecture}

\theoremstyle{definition}
\newtheorem{defn}{Definition}[section]

\newtheorem{prob}[equation]{Problem}

\newenvironment{rem}{%
\bigskip
\noindent \textsl{{\sl Remark. }}}{\bigskip}
\newenvironment{rems}{%
\bigskip
\noindent \textsl{{\sl Remarks. }}}{\bigskip}

\newcounter {own}
\def\theown {\thesection       .\arabic{own}}

\newenvironment{pf}[1][]{%
 \vskip 3mm
 \noindent
 \ifthenelse{\equal{#1}{}}%
  {{\slshape Proof. }}%
  {{\slshape #1.} }%
 }%
{\qed\bigskip}

\newcounter{alphabet}
\newcounter{tmp}
\newenvironment{Thm}[1][]{\refstepcounter{alphabet}%
\bigskip%
\noindent%
{\bf Theorem \Alph{alphabet}}%
\ifthenelse{\equal{#1}{}}{}{ (#1)}%
{\bf .} \itshape}{\vskip 8pt}

\newcommand{\IR}{{\mathbb R}}
\newcommand{\ID}{{\mathbb D}}
\newcommand{\IN}{{\mathbb N}}
\newcommand{\IC}{{\mathbb C}}

\newcommand{\CC}{{\mathcal C}}

\newcommand{\IB}{{\mathcal B}}
\newcommand{\IH}{{\mathcal H}}

\newcommand{\dist}{{\operatorname{dist}}}

\def\be{\begin{equation}}
\def\ee{\end{equation}}

\newcommand{\bee}{\begin{enumerate}}
\newcommand{\eee}{\end{enumerate}}

\newcommand{\blem}{\begin{lem}}
\newcommand{\elem}{\end{lem}}
\newcommand{\bthm}{\begin{thm}}
\newcommand{\ethm}{\end{thm}}
\newcommand{\bcor}{\begin{cor}}
\newcommand{\ecor}{\end{cor}}
\newcommand{\beg}{\begin{examp}}
\newcommand{\eeg}{\end{examp}}
\newcommand{\begs}{\begin{examples}}
\newcommand{\eegs}{\end{examples}}
\newcommand{\bdefe}{\begin{defn}}
\newcommand{\edefe}{\end{defn}}
\newcommand{\bprob}{\begin{prob}}
\newcommand{\eprob}{\end{prob}}
\newcommand{\bei}{\begin{itemize}}
\newcommand{\eei}{\end{itemize}}

\newcommand{\bcon}{\begin{conj}}
\newcommand{\econ}{\end{conj}}
\newcommand{\bcons}{\begin{conjs}}
\newcommand{\econs}{\end{conjs}}
\newcommand{\bprop}{\begin{propo}}
\newcommand{\eprop}{\end{propo}}
\newcommand{\br}{\begin{rem}}
\newcommand{\er}{\end{rem}}
\newcommand{\brs}{\begin{rems}}
\newcommand{\ers}{\end{rems}}
\newcommand{\bo}{\begin{obser}}
\newcommand{\eo}{\end{obser}}
\newcommand{\bos}{\begin{obsers}}
\newcommand{\eos}{\end{obsers}}
\newcommand{\bpf}{\begin{pf}}
\newcommand{\epf}{\end{pf}}
\newcommand{\ba}{\begin{array}}
\newcommand{\ea}{\end{array}}
\newcommand{\beq}{\begin{eqnarray}}
\newcommand{\beqq}{\begin{eqnarray*}}
\newcommand{\eeq}{\end{eqnarray}}
\newcommand{\eeqq}{\end{eqnarray*}}

\newcommand{\ov}{\overline}

\newcounter{minutes}
\divide\time by 60
\newcounter{hours}
\multiply\time by 60 \addtocounter{minutes}{-\time}
\begin{document}
\title{Bohr phenomenon for operator valued functions}
\begin{center}
{\tiny \texttt{FILE:~\jobname .tex,
        printed: \number\year-\number\month-\number\day,
        \thehours.\ifnum\theminutes<10{0}\fi\theminutes}
}
\end{center}
\author{Bappaditya Bhowmik${}^{\mathbf{*}}$}
\address{Bappaditya Bhowmik, Department of Mathematics,
Indian Institute of Technology Kharagpur, Kharagpur - 721302, India.}
\email{bappaditya@maths.iitkgp.ac.in}
\author{Nilanjan Das}
\address{Nilanjan Das, Department of Mathematics,
Indian Institute of Technology Kharagpur, Kharagpur - 721302, India.}
\email{nilanjan@iitkgp.ac.in}

\keywords{Bohr radius, harmonic function, biholomorphic function.\newline
${}^{\mathbf{*}}$ Corresponding author}

\begin{abstract}
In this article we establish Bohr inequalities for operator valued functions, which can be viewed as the analogues
of a couple of interesting results from scalar valued settings. Some results of this paper are motivated by the classical
flavor of Bohr inequality, while the others are based on a generalized concept of the Bohr radius problem.
\end{abstract}
\thanks{}

\maketitle
\pagestyle{myheadings}
\markboth{B. Bhowmik, N. Das}{Bohr phenomenon for operator valued functions}

\bigskip
\section{Introduction}
The following remarkable result was proved by Harald Bohr \cite{Bohr} in 1914.

\begin{Thm}\label{P2TheoA}
Let $f(z)= \sum_{n=0}^{\infty}a_{n}z^n$ be holomorphic in the open unit disk $\ID$ and $|f(z)|<1$ for all $z\in\ID$, then
\be\label{P2eq54}
\sum_{n=0}^{\infty}|a_n|r^n\leq 1
\ee
for all $z\in\ID$ with $|z|=r\leq1/6$.
\end{Thm}
This constant $r\leq 1/6$ was sharpened to $r\leq 1/3$ by Wiener, Riesz and Schur independently, and the inequality
$(\ref{P2eq54})$ is popularly known as \emph{Bohr inequality} nowadays.
This theorem was an outcome of the investigation on the absolute convergence problem for Dirichlet series of the form
$\sum a_n n^{-s}$, but presently it has become an independent area of research. Bohr radius problem saw a surge of
interest from many mathematicians after it found an application to the characterization problem of Banach algebras
satisfying von Neumann inequality \cite{Dix}. A part of the subsequent research in this area is directed towards extending the Bohr
phenomenon in multidimensional framework and in more abstract settings (see, for example \cite{Aiz, Aiz1, Bo, Li, Paul2}).
Bohr phenomenon is shown to have connections with local Banach space theory (cf. \cite{Def}), and is being investigated for ordinary
and vector valued Dirichlet series also (see f.i. \cite{Balu, Def1}).

We would now give a brief overview of the approaches to extend Bohr inequality in two different settings. One of them aims at
investigating Bohr radius problem from operator theoretic perspective. To be more specific,
Bohr phenomenon has been established in \cite[Theorem 2.1]{Paul} using positivity methods for operator valued holomorphic functions,
i.e. holomorphic functions from $\ID$ to $\IB(\IH)$, where $\IB(\IH)$ is the set of bounded
linear operators on a complex Hilbert space $\IH$. Suitable assumptions in terms of operator inequalities are made to replicate the
scalar valued cases. It may be mentioned here that the inequalities recorded in \cite[Theorem 2.1]{Paul} are operator valued analogues of
the classical Bohr inequality in Theorem A. In the present article we prove Bohr inequalities of similar nature for
harmonic functions from $\ID$ to $\IB(\IH)$.

Another aspect of Bohr phenomenon thrives on considering the Bohr radius problem for a holomorphic map $g$ from $\ID$ into a domain
$\Omega\subsetneq\IC$ other than $\ID$.
The key idea to accomplish that is to identify $g$ as a member of $S(f)$, $S(f)$ being the class of functions subordinate to $f$, while
$f$ is the covering map from $\ID$ onto $\Omega$ satisfying $f(0)=g(0)$. Here we clarify that for two holomorphic functions $g$
and $f$ in $\ID$, we say that $g$ is subordinate to $f$ if there exists a function $\phi$, holomorphic in $\ID$ with $\phi(0)=0$ and
$|\phi(z)|<1$, satisfying $g=f\circ\phi$. Throughout this article we denote $g$ is subordinate to $f$ by $g\prec f$. A suitable
definition for the Bohr phenomenon of $g\in S(f)$ was given in \cite{Abu} to serve the purpose stated above, which we would briefly describe here.
Let the Taylor expansions of $f$ and $g$ in a neighborhood of origin be
\be\label{P2eq51}
f(z)=\sum_{n=0}^{\infty}a_{n}z^n
\ee
and
\be\label{P2eq52}
g(z)=\sum_{k=0}^{\infty}b_{k}z^k
\ee
respectively. We will say that $S(f)$ has Bohr phenomenon if for any $g\in S(f)$, where $f$ and $g$ have the
Taylor expansions of the form $(\ref{P2eq51})$ and $(\ref{P2eq52})$ respectively in $\ID$,
there is a $r_{0}$, $0<r_{0}\leq 1$ so that
\be\label{P2eq53}
\sum_{k=1}^{\infty}|b_{k}|r^k\leq d(f(0),\partial\Omega)\,,
\ee
for $|z|=r<r_{0}$. Here $d(f(0),\partial\Omega)$ denotes the Euclidean distance between $f(0)$ and boundary of the domain $\Omega=f(\ID)$.
To see that this definition is indeed a generalization of the classical Bohr phenomenon, we observe that
whenever $\Omega=\ID$; $d(f(0), \partial\Omega)= 1-|f(0)|$, and in this case $(\ref{P2eq53})$
reduces to $(\ref{P2eq54})$.
However, to the best of our knowledge, no attempt has been made so far to obtain operator valued analogues of the
Bohr phenomenon for complex valued functions treated according to the aforesaid definition from \cite{Abu}.
Therefore, another goal of the present article is to find the same under appropriate considerations
and necessary restrictions. More precisely, we will consider a function $f$ from $\ID$ to $\IB(\IH)$, and prove Bohr inequality
when $f$ is holomorphic and satisfies certain conditions which, when restricted to the scalar valued case, coincide with the situation that $f$ maps
$\ID$ into its exterior, i.e. $\ov{\ID}^c=\{z\in\IC:|z|>1\}$. Also we prove the Bohr phenomenon for any $g\in S(f)$ when $f$ is convex or
starlike biholomorphic function. Here we clarify that, given two complex Banach spaces $X$ and $Y$ and a domain $D\subset X$, a
holomorphic mapping $f:D\to Y$ is said to be \emph{biholomorphic} on $D$ if $f(D)$ is a domain in $Y$, and $f^{-1}$ exists and is
holomorphic on $f(D)$. 
A biholomorphic function $f$ is said to be starlike on its domain $D$
with respect to $z_0\in D$ if $f(D)$ is a starlike domain with respect to $f(z_0)$, i.e.
$(1-t)f(z_0)+tf(z)\in f(D)$ for all $z\in D$ and $t\in[0,1]$ and $f$ is called \emph{starlike biholomorphic} on $D$
if $f$ is starlike with respect to $0\in D$ and $f(0)=0$. Now a biholomorphic function $f$
defined in $D$ is said to be \emph{convex} if
$f$ is starlike with respect to all $z\in D$.
In particular here we will work with $D=\ID$, $X=\IC$ and $Y=\IB(\IH)$.
It may be noted that the definition of subordination and the
class $S(f)$ for operator valued holomorphic functions can be adopted from the scalar case without any change. Now we fix some notations
for the rest of our discussions. For any $A\in\IB(\IH)$, $\|A\|$ will always denote the operator norm of $A$, and $A^*$ is the adjoint of $A$.
The operators $\mbox{Re}(A):= (A+ A^*)/2$, $\mbox{Im}(A):=(A-A^*)/2i$ and $|A|:=(A^*A)^{1/2}$
bear their usual meaning, while $B^{1/2}$ denotes the unique positive square root of a positive operator $B$. Also $\sigma(A)$ will be recognized as the
spectrum of $A$, i.e. the set of all $\lambda\in\IC$ such that $A-\lambda I$ is non invertible, $I$ being the identity operator on $\IH$.
\section{Main results}
A function $f:\ID\to\IB(\IH)$ is harmonic if and only if
\be\label{P2eq1}
f(z)=\sum_{n=0}^{\infty} A_n z^n + \sum_{n=1}^\infty B^*_n \ov{z}^n
\ee
where $A_n, B_n\in\IB(\IH)$ for all $n\in\IN\cup\{0\}$, and the series converges absolutely and locally uniformly in $\ID$
(see, for example \cite[Sec 2.4, p. 352]{Hen}). Bohr inequalities for complex valued harmonic functions have already been obtained in
\cite[Theorem 2]{Abu}. The aim of the first theorem is to derive inequalities of similar nature for operator valued harmonic functions.
Therefore we need to establish the following analogue of the Cauchy-Schwarz inequality.
\blem\label{P2lem1}
Let $\{H_n\}_{n=0}^\infty$ be a sequence in $\IB(\IH)$ such that $\sum_{n=0}^\infty |H_n|^2\in\IB(\IH)$. Then for any fixed $r\in [0,1)$,
and for any fixed non negative integer $k$,
\be\label{P2eq18}
\sum_{n=k}^\infty |H_n|r^n\leq \frac{r^k}{\sqrt{1-r^2}}\left(\sum_{n=k}^\infty|H_n|^2\right)^{1/2}.
\ee
\elem
\bpf
For any fixed $m\in\IN$ such that $m>k$, and for any $x\in\IH$, it is immediately seen that
\be\label{P2eq19}
\Big\langle\left(\sum_{n=k}^m|H_n|r^n\right)^2 x, x\Big\rangle=\left\|\sum_{n=k}^m r^n|H_n|x\right\|^2\leq \left(\sum_{n=k}^m r^n\|H_n x\|\right)^2.
\ee
for any fixed $r\in[0,1)$.
Now a use of the Cauchy-Schwarz inequality on the right hand side of inequality $(\ref{P2eq19})$ will yield
\be\label{P2eq20}
\Big\langle\left(\sum_{n=k}^m|H_n|r^n\right)^2 x, x\Big\rangle \leq \left(\sum_{n=k}^m\|H_nx\|^2\right)\left(\sum_{n=k}^m r^{2n}\right),
\ee
which would imply that
\be\label{P2eq21}
\Big\langle\left(\sum_{n=k}^m|H_n|r^n\right)^2 x, x\Big\rangle \leq \Big\langle\sum_{n=k}^m |H_n|^2x,x\Big\rangle\frac{r^{2k}}{1-r^2}.
\ee
Letting $m\to\infty$ in $(\ref{P2eq21})$ we get
\be\label{P2eq45}
\left(\sum_{n=k}^\infty |H_n|r^n\right)^2\leq \frac{r^{2k}}{1-r^2}\left(\sum_{n=k}^\infty |H_n|^2\right),
\ee
from which $(\ref{P2eq18})$ will follow (cf. \cite[p. 244, Ex. 12]{Con}).
\epf

We now state the first theorem of this article after all these preparations.
\bthm\label{P2thm1}
Let $f:\ID\to\IB(\IH)$ be a harmonic function with an expansion $(\ref{P2eq1})$ such that $\|f(z)\|\leq 1$ for each $z\in\ID$.
Then
\bee
\item[(i)]~
$|\emph{Re}(e^{i\mu}A_0)|+ \sum_{n=1}^\infty |e^{i\mu}A_n + e^{-i\mu}B_n|r^n \leq (\sqrt{1+3r^2}/\sqrt{1-r^2})I$
for $|z|=r \in [0,1)$, and for any $\mu\in\IR$.
In particular, if we assume in addition that $e^{i\mu}A_n+e^{-i\mu}B_n$ is
normal for each $n\in\IN$, then the quantity in the right hand side of the above inequality can be replaced by
$(\sqrt{1+r^2}/\sqrt{1-r^2})I$

\item[(ii)]~
$\sum_{n=1}^\infty \|e^{i\mu}A_n + e^{-i\mu}B_n\|r^n\leq \|I-\emph{Re}(e^{i\mu}A_0)\|$
for $|z|=r\leq 1/5$, and for any $\mu\in\IR$.
Moreover if we take $e^{i\mu}A_n+e^{-i\mu}B_n$ to be normal for each $n\in\IN$, then the above
inequality will hold for $r\leq 1/3$ instead of $r\leq 1/5$

\item[(iii)]~
$\sum_{n=1}^\infty |A_n|r^n+ \sum_{n=1}^\infty |B^*_n|r^n \leq (1/2)I$ for $|z|=r\leq 1/3$.
\eee
\ethm
\bpf
(i)~It is easy to observe that for each $z\in\ID$, and for
any $\mu\in\IR$,
$$
|\textrm{Re}(e^{i\mu}f(z))|^2+|\textrm{Im}(e^{i\mu}f(z))|^2=(1/2)(f(z)f(z)^*+f(z)^*f(z)).
$$
We here note that for any $A\in\IB(\IH)$, $\langle |A|^2 x,x\rangle =\|Ax\|^2 \leq\|A\|^2\langle x,x\rangle$
for any $x\in\IH$, i.e. $|A|^2\leq \|A\|^2 I$.
Using this fact, and that $\|A\|=\|A^*\|$ for any $A\in\IB(\IH)$,
we obtain $|\textrm{Re}(e^{i\mu}f(z))|^2+|\textrm{Im}(e^{i\mu}f(z))|^2\leq \|f(z)\|^2 I$,
and therefore
\be\label{P2eq2}
|\textrm{Re}(e^{i\mu}f(z))|^2\leq I.
\ee
Now
\be\label{P2eq3}
\textrm{Re}(e^{i\mu}f(z))=
\textrm{Re}(e^{i\mu}A_0)+(1/2)\sum_{n=1}^\infty(P_n z^n+ P_n^* \ov{z}^n)
\ee
where $P_n=e^{i\mu}A_n + e^{-i\mu}B_n$.
Now from $(\ref{P2eq2})$ we can write, for any $z=re^{i\theta}\in\ID$ and for any $x\in\IH$,
$$
\langle (\textrm{Re}(e^{i\mu}f(re^{i\theta})))^*(\textrm{Re}(e^{i\mu}f(re^{i\theta})))x, x\rangle\leq \langle x,x\rangle.
$$
We plug the expression $(\ref{P2eq3})$ in and fix $r\in[0,1)$ in the above inequality, and thereafter integrating both sides
of this inequality over $\theta$ from $0$ to $2\pi$ we get
$$
\langle|\textrm{Re}(e^{i\mu}A_0)|^2x,x \rangle + (1/4)\sum_{n=1}^\infty(\langle P_n^*P_nx,x\rangle + \langle P_nP_n^*x,x\rangle)r^{2n}
\leq\langle Ix,x\rangle.
$$
Therefore we conclude
$$
|\textrm{Re}(e^{i\mu}A_0)|^2+(1/4)\sum_{n=1}^\infty(|P_n|^2+|P^*_n|^2)\leq I,
$$
which implies
\be\label{P2eq6}
\sum_{n=1}^\infty|P_n|^2\leq 4(I-|\textrm{Re}(e^{i\mu}A_0)|^2).
\ee
Hence a direct use of Lemma \ref{P2lem1} (with $H_n=P_n$, $k=1$) gives
\be\label{P2eq4}
|\textrm{Re}(e^{i\mu}A_0)|+ \sum_{n=1}^\infty |e^{i\mu}A_n + e^{-i\mu}B_n|r^n \leq T+\frac{2r}{\sqrt{1-r^2}}(I-T^2)^{1/2}
\ee
where $T=|\textrm{Re}(e^{i\mu}A_0)|$. The first half of part (i) of our theorem will now follow from a computation similar to the proof
of \cite[Theorem 2.1, part 4]{Paul}, applied to $(\ref{P2eq4})$. For the sake of completion we include brief details of
the calculation. Considering the real valued function $\psi(x)= x+ (2r/\sqrt{1-r^2})\sqrt{1-x^2}$ on the interval $[0,1]$,
we see that $\psi$ attains its maximum at $x_0= \sqrt{1-r^2}/\sqrt{1+3r^2}$, and that $\psi(x)\leq \psi(x_0)= \sqrt{1+3r^2}/\sqrt{1-r^2}$
for any $x\in [0,1]$. This validates our first assertion. Further, if we assume $P_n$ is normal for each $n\in\IN$, then $|P_n|^2=|P_n^*|^2$,
which implies that the inequality $(\ref{P2eq6})$ can be improved to $\sum_{n=1}^\infty|P_n|^2\leq 2(I-|\textrm{Re}(e^{i\mu}A_0)|^2)$.
Rest of the proof can be completed by following the similar lines of computation as we did for the previous one.

(ii)~In order to establish the second part of this theorem, we first observe that if $K(z)=e^{i\mu}A_0+\sum_{n=1}^\infty P_nz^n$, then
$(\ref{P2eq3})$ implies that $\mbox{Re}(K(z))=\mbox{Re}(e^{i\mu}f(z))$, and hence from $(\ref{P2eq2})$ we get
$\|\textrm{Re}(K(z))\|\leq 1$. Now considering $\hat{K}(z)=\langle K(z)x,x\rangle$ for any fixed $x\in\IH$ with $\|x\|=1$, it is easily seen
that $|\textrm{Re}(\hat{K}(z))|\leq 1$. Therefore, $\hat{K}$ is holomorphic in $\ID$ with an expansion
$$
\hat{K}(z)=\langle e^{i\mu}A_0x,x\rangle +\sum_{n=1}^\infty\langle P_nx,x\rangle z^n
$$
which maps $\ID$ into the vertical strip $|\textrm{Re}(z)|\leq 1$. As a consequence
$$
|\langle P_nx,x\rangle|\leq 2(1-|\textrm{Re}(\langle e^{i\mu}A_0x,x\rangle)|)
$$
for all $n\in\IN$ (see, f.i. \cite[Lemma 3]{Abu}). Further, using triangle inequality we obtain
\be\label{P2eq40}
|\langle P_nx,x\rangle|\leq 2|\langle (I-\textrm{Re}(e^{i\mu}A_0))x,x\rangle|.
\ee
Taking supremum over $x\in\IH$ with $\|x\|=1$ on both sides of the inequality $(\ref{P2eq40})$,
we get
$$
\sup_{\|x\|=1}|\langle P_nx, x\rangle|\leq 2\|I-\textrm{Re}(e^{i\mu}A_0)\|,
$$
and replacing $A$ by $P_n$ in \cite[Theorem 1.2]{Gol} we have
$$
\sup_{\|x\|=1}|\langle P_nx, x\rangle|\geq (1/2)\|P_n\|.
$$
Combining the above two results we obtain
\be\label{P2eq41}
\|P_n\|\leq 4\|I-\textrm{Re}(e^{i\mu}A_0)\|
\ee
for all $n\in\IN$. The first half of part (ii) now follows from $(\ref{P2eq41})$. Now it is known that
$\|P_n\|=\sup \{|\langle P_n x,x\rangle|: x\in\IH, \|x\|=1\}$ whenever $P_n$ is normal (cf. \cite[p. 266]{Gol}, and replace $A$ by $P_n$).
Hence from $(\ref{P2eq40})$ we obtain $\|P_n\|\leq 2\|I-\textrm{Re}(e^{i\mu}A_0)\|$, which will prove the second
assertion of part (ii).

(iii)~Finally, since $\|f(z)\|\leq 1$ if and only if $|f(z)|^2\leq I$, using methods similar to the proof of part (i)
we are able to deduce
\be\label{P2eq22}
|A_0|^2+\sum_{n=1}^\infty(|A_n|^2+|B_n^*|^2)\leq I.
\ee
Observing that $|A_n|^2+|B_n^*|^2\geq (1/2)(|A_n|+|B_n^*|)^2$ for all $n\in\IN$, $(\ref{P2eq22})$ yields
\be\label{P2eq9}
\sum_{n=1}^\infty(|A_n|+|B_n^*|)^2\leq 2I.
\ee
Therefore applying Lemma \ref{P2lem1} (letting $H_n=|A_n|+|B_n^*|$, $k=1$) and $(\ref{P2eq9})$ together, we get
\be\label{P2eq10}
\sum_{n=1}^\infty |A_n|r^n+ \sum_{n=1}^\infty |B^*_n|r^n \leq \left(\frac{r\sqrt{2}}{\sqrt{1-r^2}}\right)I,
\ee
from which part (iii) will directly follow.
\epf

\brs
In connection with the above theorem the following observations are made:\\
\noindent (i)\,
Under the assumption that $e^{i\mu}A_n + e^{-i\mu}B_n\,(\mu\in\IR)$ is normal for each $n\in\IN$, from part (ii) of Theorem~\ref{P2thm1} we can write
$$
\|\textrm{Re}(e^{i\mu}A_0)\|+\sum_{n=1}^\infty \|e^{i\mu}A_n + e^{-i\mu}B_n\|r^n\leq \|\textrm{Re}(e^{i\mu}A_0)\|+\|I-\textrm{Re}(e^{i\mu}A_0)\|
$$
for $r\leq 1/3$. When restricted to scalar case, this inequality reduces to
\be\label{P2eq43}
|\textrm{Re}(e^{i\mu}A_0)|+\sum_{n=1}^\infty |e^{i\mu}A_n + e^{-i\mu}B_n|r^n\leq |\textrm{Re}(e^{i\mu}A_0)|+|1-\textrm{Re}(e^{i\mu}A_0)|
\ee
for $r\leq 1/3$, where the coefficients $A_0, A_n, B_n$ are complex numbers. Now since without loss of
generality we may consider $\textrm{Re}(e^{i\mu}A_0)\geq 0$, therefore the second part of
\cite[Theorem 2] {Abu} follows directly from $(\ref{P2eq43})$.\\
\noindent (ii)\,
Part (iii) of Theorem~\ref{P2thm1} can be thought of as an operator valued analogue of the very recent
result from \cite[p. 867, Sec. 4.4]{Kay} which improves the first part of \cite[Theorem 2]{Abu}.\\
\noindent (iii)\,
If we set $B_n=0$ for all $n\in\IN$ in $(\ref{P2eq1})$, i.e. $f$ is taken to be a holomorphic function from
$\ID$ to $\IB(\IH)$ with expansion $f(z)=\sum_{n=0}^\infty A_nz^n$, then
$(\ref{P2eq22})$ takes the form $\sum_{n=0}^\infty|A_n|^2\leq I$. Therefore an application of Lemma~$\ref{P2lem1}$
(with $H_n= A_n$, $k=0$) yields
\be\label{P2eq55}
\sum_{n=0}^\infty |A_n|r^n\leq\left(\frac{1}{\sqrt{1-r^2}}\right)I, r\in[0,1).
\ee
We observe that $1/\sqrt{1-r^2}\leq (1+ r^2/(1-r)^2)^{1/2}$, and therefore
$(\ref{P2eq55})$ is an improvement over the inequality recorded in \cite[Remark 2.2]{Paul}. Moreover from the
scalar valued results (compare \cite[Theorem 1.1]{En}), we observe that
the quantity $1/\sqrt{1-r^2}$ in inequality $(\ref{P2eq55})$ is the ``best possible", in the sense that
for the function $f(z)=((z-1/\sqrt{2})/(1-z/\sqrt{2}))I, z\in\ID$, equality occurs in $(\ref{P2eq55})$ at $r=1/\sqrt{2}$.
\ers

In the next result we establish operator valued analogue of Bohr inequality for holomorphic mappings
from $\ID$ into the exterior of $\ID$, i.e. $\ov{\ID}^c=\{z\in\IC:|z|>1\}$ (cf. \cite[Theorem 2.1]{Abu1}). In order to prove this, we
now introduce the notions of the spherical and the Hausdorff distance. Let $\hat{\IC}=\IC\cup\{\infty\}$ be the extended complex plane.
The spherical distance $\lambda$ between two points $z_1, z_2\in\hat{\IC}$ is given by
$$
\lambda(z_1,z_2)=
\begin{cases}
\frac{|z_1-z_2|}{\sqrt{1+|z_1|^2}\sqrt{1+|z_2|^2}},\, \mbox{if}\, z_1, z_2\in\IC,\\
\frac{1}{\sqrt{1+|z_1|^2}},\, \mbox{if}\, z_2=\infty.\\
\end{cases}
$$
Also it is well known that the collection $\CC$ of compact subsets of $\IC$ is a metric space with respect to the Hausdorff distance
$d_h$ given by
$$
d_h(A, B)= \max\{\sup_{x\in A}\dist(x, B),\,\sup_{x\in B}\dist(x, A)\}, \, A, B\in \CC,
$$
where $\dist(p, E):=\inf\{|p-e|:e\in E\}$ for any $E\subset\IC$ and for any $p\in\IC$. Now since for any $A\in\IB(\IH)$, $\sigma(A)\in \CC$, we
are able to consider the mapping $A\mapsto\sigma(A)$ from $\IB(\IH)$ to the metric space$(\CC, d_h)$, which is continuous on the subset of
normal operators, equipped with the operator norm (see f.i. \cite{New}).
\bthm\label{P2thm2}
Suppose $f:\ID\to\IB(\IH)$ be holomorphic with an expansion
\be\label{P2eq23}
f(z)=\sum_{n=0}^\infty A_n z^n, z\in\ID
\ee
such that $|f(z)|>I$ for all $z\in\ID$. Also suppose $f(z)$ is
normal for each $z\in\ID$, $f(0)=A_0> 0$ and $\sigma(f(z))$ does not separate $0$ from $\infty$ for any $z\in\ID$. Then
\be\label{P2eq11}
\lambda\left(\sum_{n=0}^\infty\|A_n\|r^n, \|A_0\|\right)\leq \lambda\left(\|A_0\|,1\right)
\ee
for $|z|=r\leq (2(\log\|A_0\|/\|\log A_0\|)-1)/(2(\log\|A_0\|/\|\log A_0\|)+1)$.
\ethm
\bpf
Since $|f(z)|>I$, we have $\langle |f(z)|x, x\rangle > \langle x, x\rangle$ for any $x\in\IH\setminus\{0\}$, and for each $z\in\ID$. A use of
the Cauchy-Schwarz inequality exhibits that $\|f(z)x\|>\|x\|$, which further implies that $\|(f(z)-\lambda I)x\|>(1-|\lambda|)\|x\|$
for any $\lambda\in\IC$, i.e. $f(z)-\lambda I$ is bounded below for any $\lambda\in\ID$. As $f(z)$ is normal,
$\sigma(f(z))\subset \ID^c$ for each $z\in\ID$. Since $\sigma(f(z))$ does not separate $0$ from $\infty$, therefore it is possible
to choose a holomorphic single valued branch of complex logarithm on a simply connected domain
$\Delta_z$ that contains $\sigma(f(z))$, but does not contain 0.
As a consequence we are able to define $\log f(z)$ as follows:
\be\label{P2eq50}
\log f(z)=\frac{1}{2\pi i}\int_{\Gamma}(\log \xi)(\xi I-f(z))^{-1} d\xi, z\in\ID
\ee
where $\Gamma$ is a system of closed, positively oriented, rectifiable curves inside $\Delta_z$ which encloses $\sigma(f(z))$
(cf. \cite[pp. 199-201]{Con}). Now it is also known that for each fixed $z\in\ID$,
$\log f(z)$ is normal, and $(\log f(z))^*=F(f(z)^*)$, where $F(z)=\ov{\log \ov{z}}$ (see f.i. \cite[p. 205, Ex. 7, 8]{Con}).
As $\exp{z}$ is an entire function and $\exp(\ov{\log \ov{z}})=z$, it follows that $\exp((\log f(z))^*)=f(z)^*$ (see \cite[p. 205, Ex. 4]{Con}).
As a consequence of these facts, we obtain $\exp(2 \mbox{Re}(\log f(z)))= f(z)^*f(z)$. It is easy to see that for any $x\in\IH\setminus\{0\}$,
$$
\langle\exp(2 \mbox{Re}(\log f(z))) x,x\rangle=\langle f(z)^*f(z)x, x\rangle=\|f(z)x\|^2 >\|x\|^2,
$$
which, after an application of the Cauchy-Schwarz inequality asserts that
$$\|\exp(2 \mbox{Re}(\log f(z)))x\|> \|x\|.
$$
Therefore $\sigma(\exp(2 \mbox{Re}(\log f(z)))x)\subset\ID^c$, and since the
operator $\exp(2 \mbox{Re}(\log f(z)))$ is positive, we conclude that
$\sigma(\exp(2 \mbox{Re}(\log f(z))))\subset [1, \infty)$.
Now we know that $\sigma(2\mbox{Re}(\log f(z)))\subset \IR$, and hence $\exp(\sigma(2\mbox{Re}(\log f(z))))\subset (0, \infty)$.
As a result, choosing the principal branch of complex logarithm over the slit plane $\IC\setminus (-\infty, 0]$,
we get $\log(\exp(2 \mbox{Re}(\log f(z))))= 2 \mbox{Re}(\log f(z))$.
Now applying the spectral mapping theorem, we conclude that
$$
\sigma(2 \mbox{Re}(\log f(z)))=\log\left(\sigma(\exp(2 \mbox{Re}(\log f(z))))\right)\subset [0, \infty).
$$
As $2 \mbox{Re}(\log f(z))$ is self adjoint, $2 \mbox{Re}(\log f(z))\geq 0$. Moreover, as $A_0>0$,
$\sigma(A_0)\subset [1, \infty)$. Hence to define $\log A_0$ from $(\ref{P2eq50})$, we choose, in particular, the principal branch of complex logarithm
on the simply connected domain $\Delta_0=\IC\setminus (-\infty, 0 ]$ containing $\sigma(A_0^*)=\sigma(A_0)$. Now as
$F(z)= \log z$ over $[1, \infty)$, $F(z)= \log z$,
$z\in\IC\setminus (-\infty, 0 ]$.
Therefore, $(\log A_0)^*= \log A^*_0= \log A_0$, which in turn gives
$\log A_0\geq 0$.
Our aim is now to show that $\log f(z)$ is holomorphic at each $z\in\ID$.
As $f(z)$ is holomorphic, and therefore continuous on $\ID$, $\lim_{h\to 0}\|f(z+h)-f(z)\|=0$. Since $f(z)$ is also normal for each $z\in\ID$,
we have
$$
\lim_{h\to 0}d_h(\sigma(f(z+h)), \sigma(f(z)))=0,
$$
thus we infer that for any $h\in\IC$ with $|h|$ small enough, $\sigma(f(z+h))$ is enclosed by $\Gamma$ again. As a result we are able to
show that the limit
$$
\lim_{h\to 0} \frac{1}{2\pi ih}\int_\Gamma (\log \xi)((\xi I-f(z+h))^{-1}-(\xi I-f(z))^{-1}) d\xi
$$
exists and is equal to
$$
\frac{1}{2\pi i}\int_{\Gamma}(\log \xi)(\xi I-f(z))^{-1} f^\prime(z) (\xi I-f(z))^{-1} d\xi ,
$$
thereby proving that $\log f(z)$ is holomorphic in $\ID$.
In view of the above discussion, there exist a Hilbert space $\mathcal{K}$,
a unitary operator $U$ on $\mathcal{K}$ and a bounded linear operator $V:\IH\to\mathcal{K}$ such that
$2\log A_0=V^*V$ and $C_n=V^*U^nV$ for all $n\geq 1$, where $\log f(z)= \log A_0+\sum_{n=1}^\infty C_nz^n$, $z\in\ID$
(see f.i. \cite[Ex. 3.15, 3.16, 4.14]{Paul1}). Hence for any $z\in\ID$, we have
$$
2\log f(z)=V^*(I+zU)(I-zU)^{-1}V,
$$
which immediately gives
\be\label{P2eq14}
f(z)=\exp\left((1/2)V^*(I+zU)(I-zU)^{-1}V\right).
\ee
From $(\ref{P2eq14})$ it can be observed that all the $A_n$ 's are the combinations of $U, V$ and $V^*$, associated with
nonnegative real constants only. Therefore a use of triangle inequality will provide the upper bounds for $\|A_n\|$ 's,
which are the combinations of $\|U\|=1$, $\|V\|=\|V^*\|$, associated with the same constants. Hence after appropriate
rearrangement we find that for any $|z|=r$,
\be\label{P2eq15}
\sum_{n=0}^\infty\|A_n\|r^n\leq \sum_{n=0}^\infty \frac{1}{n!} \left(\frac{\|V\|^2}{2}\frac{1+r}{1-r}\right)^n=
\exp\left(\frac{\|V\|^2}{2}\frac{1+r}{1-r}\right).
\ee
As $\|V\|^2=2\|\log A_0\|$, therefore we get
\be\label{P2eq16}
\sum_{n=0}^\infty\|A_n\|r^n\leq \exp(2\log \|A_0\|)=\|A_0\|^2
\ee
whenever $r\leq r_0:=(2(\log\|A_0\|/\|\log A_0\|)-1)/(2(\log\|A_0\|/\|\log A_0\|)+1).$
Now if $\alpha, \beta, \gamma$ are nonnegative real numbers satisfying
$\gamma\leq \alpha\leq \beta$, then it is easily seen that
$$
(\alpha-\gamma)^2(1+\beta^2)-(\beta-\gamma)^2(1+\alpha^2)= (\alpha-\beta)((\alpha-\gamma)+ (\beta-\gamma)+
\alpha\gamma(\beta-\gamma)+ \beta\gamma(\alpha-\gamma))\leq 0.
$$
As a consequence,
$(\alpha-\gamma)/\sqrt{1+\alpha^2}\leq (\beta-\gamma)/\sqrt{1+\beta^2}$,
which readily gives
\be\label{P2eq17}
\lambda(\alpha, \gamma)\leq\lambda(\beta, \gamma).
\ee
Setting $\alpha= \sum_{n=0}^\infty\|A_n\|r^n$, $\beta= \|A_0\|^2$ and $\gamma=\|A_0\|$,
we observe that $\gamma\leq\alpha\leq\beta$ if $r\leq r_0$, and therefore from $(\ref{P2eq17})$ we get
$$
\lambda\left(\sum_{n=0}^\infty\|A_n\|r^n, \|A_0\|\right)\leq \lambda\left(\|A_0\|^2, \|A_0\|\right)
$$
for $r\leq r_0$. A little computation using the AM-GM inequality yields
\beqq
\lambda(\|A_0\|^2,\|A_0\|)&\leq &\|A_0\|(\|A_0\|-1)/(\sqrt{1+\|A_0\|^2}\sqrt{2}\|A_0\|)\\
&=&(\|A_0\|-1)/(\sqrt{2}\sqrt{1+\|A_0\|^2})=\lambda(\|A_0\|,1).
\eeqq
It is now clear that an application of the above inequality upon the right hand side of the previous one
will complete the proof.
\epf

\br
It does not seem plausible that we can get a uniform bound on $|z|$ which is not dependent on $A_0$ and
will still imply $(\ref{P2eq11})$. Nevertheless, if $f$ is taken to be scalar
valued, then since it is always possible to assume that $f(0)> 0$, the quantity
$(2(\log\|A_0\|/\|\log A_0\|)-1)/(2(\log\|A_0\|/\|\log A_0\|)+1)$ converts to the constant $1/3$, and
$\lambda(\|A_0\|,1)=\lambda(A_0,\partial\Omega)$, $A_0$ being an element of $\IC$ and $\partial\Omega$ being the boundary of $\ov{\ID}^c$.
Therefore Theorem~\ref{P2thm2} provides an operator valued analogue of \cite[Theorem 2.1]{Abu1}. It is interesting to note that here one has
to consider spherical distance between complex numbers to obtain Bohr inequality instead of the Euclidean distance used in
$(\ref{P2eq53})$.
\er

We will now discuss the operator valued analogues of Bohr radius problem for the subordination classes of functions
which belong to well known subclasses of scalar valued univalent functions. We therefore consider $f$ to be biholomorphic for our purpose.
Now it is possible to carry out further investigation if we restrict $f$ to some subclass of biholomorphic functions.
In particular we intend to establish Bohr inequalities for $g\in S(f)$ where $f:\ID\to\IB(\IH)$
is a convex or starlike biholomorphic function.
Apart from the definitions given in the introduction,
the reader is urged to glance through \cite{Kohr} for a rich exposition of Banach space valued starlike and convex
biholomorphic functions.
For our purpose we suppose that $g\in S(f)$
has an expansion
\be\label{P2eq24}
g(z)=\sum_{k=0}^\infty B_kz^k, z\in\ID.
\ee
Also we mention that for any scalar valued univalent function $F$ defined on $\ID$,
the Euclidean distance between $F(0)$ and the boundary $\partial\Omega$ of $\Omega=F(\ID)$ is given
by $d(F(0), \partial\Omega)=\lim\inf_{|z|\to1-}|F(z)-F(0)|$,
which will be used frequently in our forthcoming discussions.
\bthm\label{P2thm3}
Let $f:\ID\to\IB(\IH)$ be a convex biholomorphic function and $g\in S(f)$ with expansions $(\ref{P2eq23})$
and $(\ref{P2eq24})$ respectively. Then for $|z|=r\leq 1/(1+2\|A_1\|\|A_1^{-1}\|)$ we have
\be\label{P2eq46}
\sum_{k=1}^\infty \|B_k\|r^k\leq \lim\inf_{|z|\to1-}\|f(z)-f(0)\|.
\ee
Also for $|z|=r\leq 1/3$ we have
\be\label{P2eq47}
\sum_{k=1}^\infty |B_k|r^k\leq (1/2)|A_1|.
\ee
\ethm
\bpf
We observe that the well known argument used in proving \cite[Theorem X]{Rogo} can be
used in a similar fashion for $g\in S(f)$ where $f$ is operator valued convex biholomorphic function.
Thus we have $B_k=\phi^\prime(0)f^\prime(0)$, $k\geq 1$
for some holomorphic map $\phi:\ID\to\ID$ with $\phi(0)=0$.
Therefore we immediately see $\|B_k\|\leq \|A_1\|$ and hence
the following inequality will hold:
\be\label{P2eq29}
\sum_{k=1}^\infty \|B_k\|r^k\leq (r/1-r)\|A_1\|.
\ee
Now for any fixed $a\in\ID$, we construct the familiar Koebe transform as follows:
\be\label{P2eq30}
G(z)=(1-|a|^2)^{-1}(f^\prime(a))^{-1}\left(f\left((z+a)(1+\ov{a}z)^{-1}\right)-f(a)\right), z\in\ID.
\ee
We see that $G(z)$ is convex biholomorphic with the normalization $G(0)=0$ and $G^\prime(0)=I$.
From \cite[Theorem 6.3.5]{Kohr} we get that $G$ satisfies
$$
zG^{\prime\prime}(z)+ G^\prime(z)=p(z)G^\prime(z),
$$
where $p:\ID\to\IC$
be holomorphic with $\textrm{Re}(p(z))>0$ for all $z\in\ID$ and $p(0)=1$. Therefore for any fixed $x\in\IH$  with $\|x\|=1$,
the function $\hat{G}:\ID\to\IC$ defined by
$$\hat{G}(z)=\langle G(z)x,x\rangle
$$
satisfies $\hat{G}(0)=\hat{G^\prime}(0)-1=0$ and $z\hat{G}^{\prime\prime}(z)+ \hat{G}^\prime(z)=p(z)\hat{G}^\prime(z)$,
which together implies $\hat{G}(z)$ is a complex valued normalized convex univalent function (see \cite[Theorem 2.2.3]{Kohr}).
As a consequence, $\lim\inf_{|z|\to1-}|\hat{G}(z)|\geq 1/2$ (cf. \cite[Theorem 2.2.9]{Kohr}), which, after an application of the Cauchy-Schwarz
inequality for inner product yields
\be\label{P2eq31}
\lim\inf_{|z|\to 1-}\left\|(1-|a|^2)^{-1}(f^\prime(a))^{-1}\left(f\left((z+a)(1+\ov{a}z)^{-1}\right)-f(a)\right)\right\|\geq 1/2.
\ee
Now inequality $(\ref{P2eq31})$ will further give
\be\label{P2eq32}
\lim\inf_{|z|\to1-}\left\|\left(f\left((z+a)(1+\ov{a}z)^{-1}\right)-f(a)\right)\right\|\geq
(1-|a|^2)/(2\|(f^\prime(a))^{-1}\|).
\ee
In particular for $a=0$,
\be\label{P2eq34}
\lim\inf_{|z|\to1-}\|f(z)-f(0)\|\geq 1/2\|A_1^{-1}\|.
\ee
From $(\ref{P2eq29})$ and $(\ref{P2eq34})$, a little computation reveals that
$(\ref{P2eq46})$ will hold if
$$
(r/1-r)\|A_1\|\leq 1/2\|A_1^{-1}\|\,~ \mbox{or equivalently if}\,~ r\leq 1/(1+2\|A_1\|\|A_1^{-1}\|).
$$
Now going back to the relation $B_k=\phi^\prime(0)f^\prime(0)$, it is readily seen that $|B_k|\leq|A_1|$ for any $k\geq 1$,
and therefore
\be\label{P2eq48}
\sum_{k=1}^\infty |B_k|r^k\leq (r/1-r)|A_1|.
\ee
It is easy to see that for $r\leq 1/3$, $(\ref{P2eq48})$ is converted to $(\ref{P2eq47})$.
\epf

\brs
We make the following observations related to Theorem~\ref{P2thm3}.

\noindent (i)\, The quantity $1/(1+2\|A_1\|\|A_1^{-1}\|)$ in Theorem $\ref{P2thm3}$ will turn into $1/3$ for scalar valued functions, as whenever
$A_1$ is a scalar, $\|A_1\|\|A_1^{-1}\|=1$. Therefore $(\ref{P2eq46})$ gives operator valued analogue
of Bohr phenomenon for the subordinating family of a complex valued convex univalent function defined on $\ID$ (compare \cite[Remark 1]{Abu}).
\noindent (ii)\,
The right hand side of the inequality $(\ref{P2eq47})$ can be further estimated to observe $(1/2)|A_1|\leq d(f(0), \partial\Omega)$
when scalar valued functions are being considered (see \cite[Lemma 3]{Abu}), $\partial\Omega$ being the boundary of $\Omega=f(\ID)$.
Due to this fact it can be thought of as a generalization of the Bohr phenomenon mentioned in \cite[Remark 1]{Abu}.

\ers

Before we proceed further we prove the following lemma which will be required to establish the subsequent results.
\blem\label{P2lem2}
Let $f:\ID\to\IB(\IH)$ be holomorphic and $g\in S(f)$ with expansions
$(\ref{P2eq23})$ and $(\ref{P2eq24})$ respectively.
Then for $|z|=r\leq 1/3$ we have
\bee
\item[(i)]~
$\sum_{k=1}^\infty |B_k|r^k\leq (\sum_{n=1}^\infty \|A_n\|r^n)I$
\item[(ii)]
$\sum_{k=1}^\infty \|B_k\|r^k\leq \sum_{n=1}^\infty \|A_n\|r^n$.
\eee
\elem
\bpf
Since $g\prec f$, there exists a function $\phi$, holomorphic in $\ID$, satisfying $\phi(0)=0$ and $\phi(\ID)\subset\ID$ such that
\be\label{P2eq56}
g=f\circ\phi\,.
\ee
Since $\phi$ is holomorphic, the Taylor expansion of the $t$-th power of $\phi$, where $t\in\IN$, can be written as
\be\label{P2eq57}
\phi^t(z)= \sum_{l=t}^{\infty}{\alpha}{_l}^{(t)}z^l.
\ee
Now we plug equality $(\ref{P2eq57})$ into $(\ref{P2eq56})$, and equating the coefficients for $z^k$ from both sides we have, for any $k\geq 1$:
$$
B_{k}=\sum_{n=1}^k\alpha_{k}^{(n)}A_n.
$$
Now we see that
\begin{align*}
\sum_{k=1}^{m}|B_{k}|r^k &=\sum_{k=1}^{m}\left|\sum_{n=1}^k{\alpha_{k}^{(n)}A_n}\right|r^k\\
&\leq\left(\sum_{k=1}^{m}\left\|\sum_{n=1}^k{\alpha_{k}^{(n)}A_n}\right\|r^k\right)I
\leq\left(\sum_{k=1}^{m}\sum_{n=1}^k{|\alpha_{k}^{(n)}|\|A_n\|}r^k\right)I.
\end{align*}
We observe that the rightmost term of the above inequality can be written as $\left(\sum_{n=1}^m\|A_n\|M_m^{(n)}(r)\right)I$
where $M_m^{(n)}(r):=\sum_{k=n}^m|\alpha_k^{(n)}|r^k.$
The proof of part (i) can now be completed by adopting the techniques similar to the proof of \cite[Lemma 1]{BB} hereafter.
Further, part (ii) can be proved by directly following the same line of computations as in the proof of \cite[Lemma 1]{BB}.
\epf

We now state and prove a theorem including Bohr phenomenon for $S(f)$ where $f$ is an operator valued normalized
starlike biholomorphic function. It may be mentioned that the known techniques to find out the coefficient
bounds for functions subordinate to a complex valued normalized starlike univalent function do not seem to
be directly applicable in this situation, while a use of Lemma~\ref{P2lem2} will prove the following theorem.
\bthm\label{P2thm4}
Let $f:\ID\to\IB(\IH)$ be a normalized starlike biholomorphic function with an expansion $f(z)=zI+\sum_{n=2}^\infty A_nz^n$
and $g\in S(f)$ with an expansion $(\ref{P2eq24})$. Then for $|z|=r\leq 3-2\sqrt{2}$ we have
\bee
\item[(i)]~
$\sum_{k=1}^\infty \|B_k\|r^k\leq \lim\inf_{|z|\to1-}\|f(z)\|$
\item[(ii)]~
$\sum_{k=1}^\infty |B_k|r^k\leq (1/4)I$.
\eee
\ethm
\bpf
From \cite[Theorem 6.2.6]{Kohr}, it is seen that a starlike biholomorphic function $f:\ID\to\IB(\IH)$ normalized by
$f(0)=f^\prime(0)-I=0$ satisfies
\be\label{P2eq49}
zf^\prime(z)=p(z)f(z), z\in\ID
\ee
where $p:\ID\to\IC$ be holomorphic with $\textrm{Re}(p(z))>0$ for all $z\in\ID$ and $p(0)=1$.
Here we mention that a holomorphic function $f:\ID\to\IC$, normalized by $f(0)=f^\prime(0)-1=0$ is starlike univalent
if and only if $(\ref{P2eq49})$ holds.
Now a standard method based on induction
(see f.i. \cite[Theorem 2.2.16]{Kohr}) yields $\|A_n\|\leq n$ for all $n\geq 2$.
As a consequence $\sum_{n=1}^\infty\|A_n\|r^n\leq r/(1-r)^2$, where $A_1=I$. Now
let us define $G:\ID\to\IC$ by
$$
G(z)=\langle f(z)x, x\rangle
$$
where $x\in\IH$ with $\|x\|=1$.
It is easy to see that $G(0)=G^\prime(0)-1=0$.
Therefore following the similar lines of argument as in the proof of
Theorem \ref{P2thm3}, $(\ref{P2eq49})$ implies that $G$ is a starlike univalent function.
Hence $\lim\inf_{|z|\to 1-}|G(z)|\geq 1/4$ (see \cite[Theorem 1.1.5]{Kohr},
and observe that the Koebe function $k(z)=z/(1-z)^2$
which skips the value $-1/4$ is starlike univalent), and as a result the Cauchy-Schwarz inequality for inner product gives
$\lim\inf_{|z|\to 1-}\|f(z)\|\geq 1/4$.
From a direct calculation,
$\sum_{n=1}^\infty\|A_n\|r^n\leq 1/4$ for $|z|=r\leq 3-2\sqrt{2}$, which is
less than $1/3$. By virtue of the Lemma $\ref{P2lem2}$, our proofs for both part (i) and (ii) will be complete.
\epf

\brs
We end the article with the following observations:\\
(i)\, It is immediately seen that for complex valued function $f$, part (i) of the Theorem~\ref{P2thm4} converts to the Bohr inequality
for $S(f)$ where $f$ is a normalized starlike univalent function. Again, if $f$ is a complex valued normalized starlike univalent function defined on $\ID$,
the right hand side of the inequality in part (ii) is converted to $1/4$ which is known to be less or equal to $d(f(0), \partial\Omega)$, $\partial\Omega$
being the boundary of $\Omega=f(\ID)$, and thereby showing that part (ii) can also be considered as an operator valued analogue of the Bohr phenomenon for
$S(f)$. We note that the scalar valued result is a direct consequence of
\cite[Theorem 1]{Abu}.

\noindent (ii)\, In view of Theorem~\ref{P2thm3} and Theorem~\ref{P2thm4}, it is a natural question to ask
if the inequality $(\ref{P2eq46})$ holds for $|z|=r\leq r_0$ for some $r_0>0$,
where $f$ is any function in the entire family of biholomorphic functions from $\ID$ to $\IB(\IH)$ and $g\in S(f)$.
The Bohr radius $1/(1+2\|A_1\|\|A_1^{-1}\|)$
determined in the first part of the Theorem~\ref{P2thm3}
is not bounded below by a positive constant if we allow $A_1$ to be any invertible operator from $\IB(\IH)$,
$\IH$ varying on the family of complex Hilbert spaces. Therefore
we remark that the answer of the aforesaid question could possibly be negative, even for $f$ being convex biholomorphic,
and that this can be an interesting problem for future research.
However, similar problem for Banach space valued holomorphic functions in $\ID$ has already been settled (cf. \cite[Theorem 1.2]{Bl}),
where the notion of the Bohr inequality is analogous to $(\ref{P2eq54})$.
\ers

\end{document}